\documentclass[10pt,a4paper] {article}  
	
	\usepackage[english]{babel}
	
	\usepackage[margin=2cm,font=footnotesize,labelfont=bf]{caption} 
	
	\usepackage{amsmath} 
	\usepackage{amsthm} 
	\usepackage{amssymb}
	\usepackage{amsfonts}
	\usepackage[dvips]{graphicx}
	
	\usepackage[table]{xcolor}
	\usepackage{epstopdf}

	\usepackage{url}
	
	\usepackage{multirow}
	
	\usepackage{authblk}

	\DeclareGraphicsExtensions{.jpg,.pdf,.png,.gif,.eps}

	\usepackage[left=3cm,top=3cm,bottom=3cm,right=3cm]{geometry}
	
	\makeindex
	
	\title{MEIGO: an open-source software suite based on metaheuristics for global optimization in systems biology and bioinformatics}

	\author[1]{Jose A. Egea\thanks{josea.egea@upct.es}}
	\author[2]{David Henriques\thanks{davidh@iim.csic.es}}
	\author[3]{Thomas Cokelaer\thanks{cokelaer@ebi.ac.uk}}
	\author[2]{Alejandro F. Villaverde\thanks{afvillaverde@iim.csic.es}}
	\author[2]{Julio R. Banga\thanks{julio@iim.csic.es}}
	\author[3]{Julio Saez-Rodriguez\thanks{saezrodriguez@ebi.ac.uk}}
	
	\affil[1]{Department of Applied Mathematics and Statistics, Universidad Polit\'{e}cnica de Cartagena, 30202 Cartagena, Spain}
	\affil[2]{(Bio)Process Engineering Group, Spanish National Research Council, IIM-CSIC, 36208 Vigo, Spain}
	\affil[3]{European Bioinformatics Institute (EMBL-EBI), Wellcome Trust Genome Campus, Cambridge CB10 1SD, UK}

	\date{\today}

\begin{document}
	
	\maketitle
	
	\begin{abstract}
	Optimization  is key to solve many problems in computational biology. Global optimization methods provide a robust methodology, and metaheuristics in particular have proven to be the most efficient methods for many applications. Despite their utility, there is limited availability of metaheuristic tools. We present MEIGO, an R and Matlab optimization toolbox (also available in Python via a wrapper of the R version), that implements metaheuristics capable of solving diverse problems arising in systems biology and bioinformatics: enhanced scatter search method (eSS) for continuous nonlinear programming (cNLP) and mixed-integer programming (MINLP) problems, and variable neighborhood search (VNS) for Integer Programming (IP) problems. Both methods can be run on a single-thread or in parallel using a cooperative strategy. The code is supplied under GPLv3 and is available at \url{http://www.iim.csic.es/~gingproc/meigo.html}. Documentation and examples are included. The R package has been submitted to Bioconductor. We evaluate MEIGO against optimization benchmarks, and illustrate its applicability to a series of case studies in bioinformatics and systems biology, outperforming other state-of-the-art methods. MEIGO provides a free, open-source platform for optimization, that can be applied to multiple domains of systems biology and bioinformatics. It includes efficient state of the art metaheuristics, and its open and modular structure allows the addition of further methods.
	\end{abstract}

	\section*{Background}
	Mathematical optimization plays a key role in systematic decision making processes, and is used virtually in all areas of science and technology where problems can be stated as finding the best solution among a set of feasible ones. In bioinformatics and systems biology, there has been a plethora of successful applications of optimization during the last two decades (see reviews in \cite{Greenberg_Hart_Lancia:04,Larranaga_et_al:06,Festa:07, Banga:08,Sun_Garibaldi_Hodgman:12}).
	Many problems in computational biology can be formulated as IP problems, such as e.g. sequence alignment, genome rearrangement and protein structure prediction problem \cite{Greenberg_Hart_Lancia:04,Festa:07}, or design of synthetic biological networks \cite{Marchisio_Stelling:09}. Deterministic and stochastic/heuristic methods for optimization arising from the area of machine learning have also been extensively applied \cite{Larranaga_et_al:06}. In addition to combinatorial optimization, other important classes of optimization problems that have been extensively considered, especially in systems biology, are cNLP and mixed-integer dynamic optimization , targeting problems such as parameter estimation and optimal experimental design \cite{Banga_Balsa-Canto:08,Sun_Garibaldi_Hodgman:12}.
	
	A number of authors have stressed the need of using suitable global optimization methods due to the non-convex (multimodal) nature of many of these problems \cite{Moles_Mendes_Banga:03,Banga:08,Ashyraliyev_et_al:08}. Roughly speaking, global optimization methods can be classified into exact and stochastic approaches. Exact methods can guarantee convergence to global optimality, but usually the associated computational effort is prohibitive for realistic applications. In contrast, stochastic methods are often able to locate the vicinity of the global solution in reasonable computation times, but without guarantees.  Metaheuristics (i.e. guided heuristics) are a particular class of stochastic methods that have been shown to perform very well in a broad range of applications \cite{Sun_Garibaldi_Hodgman:12}.
	
	Motivated by this, we developed the software suite MEIGO (MEtaheuristics for systems biology and bIoinformatics Global Optimization) which provides state of the art metaheuristics (eSS and VNS) in open-source R and Matlab versions (also available in Python via a wrapper for the R version) covering the most important classes of problems, namely (i) problems with real-valued (cNLP’s) and mixed-integer decision variables (MINLP’s), and (ii) problems with integer and binary decision variables (IP’s).  Furthermore, MEIGO allows the user to carry out parallel computation using cooperative strategies \cite{Villaverde_et_al:12}. 
	MEIGO can optimize arbitrary objective functions that are handled as black-boxes. Thus, it is possible to apply these methods to complex problems where the objective function and/or the constraints need to solve additional problems. For example, CellNOpt \cite{Terfve_et_al:12}, SBToolbox \cite{Schmidt_Jirstrand:06}, AMIGO \cite{Balsa-Canto_Banga_11} and Potterswheel \cite{Maiwald_Eberhardt_Blumberg:12} use eSS for dynamic model calibration. Some recent successful applications of eSS in the fields of systems biology can be found in  \cite{Balsa-Canto_Alonso_Banga:10,Skanda_Lebiedz:10,Yuraszeck_et_al:10,Jia_Stephanopoulos_Gunawan:11,Heldt_Frensing_Reichl:12,Higuera_et_al:12,Jia_Stephanopoulos_Gunawan:12,MacNamara_et_al:12,Sriram_Rodriguez-Fernandez_DoyleIII:12a,Sriram_Rodriguez-Fernandez_DoyleIII:12b,Freund_et_al:13,Francis_Garcia_Middleton:press}. Besides, it was shown that eSS outperformed the various optimization methods available in the Systems Biology Toolbox \cite{Egea_Schmidt_Banga:08}.

	\section*{Methods}
	\subsection*{Enhanced Scatter Search (eSS)}
	Scatter search \cite{Glover_Laguna_Marti:00} is a population-based metaheuristics which can be classified as an evolutionary optimization method. In contrast with other popular population-based metaheuristics like e.g., genetic algorithms, the population size in scatter search is small, and the combinations among its members are carried out in a systematic way rather than randomly. The current population in scatter search is commonly named ``Reference Set'' (RefSet) and the so-called \textit{improvement method}, which consists in a local search to increase the convergence to optimal solutions, can be applied with more or less frequency to its members. A set of improvements have been implemented in the enhanced scatter search method. Among the most remarkable changes, we can mention the replacement method. Unlike in the original scatter search scheme, which uses a $\mu+\lambda$ replacement, the enhanced scatter search uses a $1+1$ replacement similar to the strategy used in a very efficient evolutionary method, Differential Evolution \cite{Storn_Price:97}. Besides, the ``go-beyond'' strategy to exploit combinations which explore promising directions has been implemented. The use of memory is also exploited to select the most efficient initial points to perform local searches, to avoid premature convergence and to perturb solution vectors which are stuck in stationary points. More details about the enhanced scatter search scheme can be found in \cite{Egea_Marti_Banga_10}.
	
	\subsection*{Variable Neighbourhood Search  (VNS)}
	Variable Neighbourhood Search is a trajectory-based metaheuristics for global optimization. It was introduced by Mladenovi\'{c} and Hansen \cite{Mladenovic_Hansen:97} and has gained popularity in recent years in the field of global optimization.
	VNS performs a local search, by evaluating the objective function around an incumbent solution, and repeats the procedure visiting different neighbourhoods to locate different local optima, among which the global optimum is expected to be found. One of the key points of the algorithm is the strategy followed to change the current neighbourhood. VNS usually seeks a new neighbourhood by perturbing a set of decision variables using a distance
	criterion. Once a new solution has been created in the new neighbourhood, a new local search is performed. The typical scheme consists of visiting neighbourhoods close to the current one (i.e., perturbing small set of solutions), until no further improvement is achieved. Then, ore distant neighbourhoods are explored. 
	Apart from this basic scheme, we have implemented advanced strategies to
	avoid cycles in the search (e.g., not repeating the perturbed decision variables in
	consecutive neighbourhood searches), to increase the efficiency when dealing with
	large-scale problems (e.g., by allowing a maximum number of perturbed decision
	variables, like in the Variable Neighbourhood Decomposition Search strategy \cite{Hansen_Mladenovic_Perez-Brito:01}), or to modify the search aggressiveness to locate high quality solutions (even if not the global optimum) in short computation times if required. Other heuristics, like the ``go-beyond'' strategy, to exploit promising directions during the local search have been adapted from other metaheuristics for continuous optimization \cite{Egea_Marti_Banga_10}.

	\subsection*{Cooperation}
	The cooperation scheme implemented in MEIGO is based in the following idea: to run, in parallel, several implementations or threads of an optimization algorithm–which may have different
	settings and/or random initializations–and exchange information between them. Since the nature of the optimization algorithms implemented in MEIGO is essentially different, we distinguish between eSS (the population based method) and VNS (the trajectory based method), following the classification
	proposed in \cite{Toulouse_Crainic_Sanso:04}:
	
	\begin{enumerate}
	\item Information available for sharing: the best solution
	found and, optionally for eSS, the \textit{RefSet}, which contains information about the diversity of
	solutions.
	\item Threads that share information: all of them.
	\item Frequency of information sharing: the threads exchange information at a fixed interval $\tau$.
	\item Number of concurrent programs: $\eta$.
	\end{enumerate}
	
	Each of the  $\eta$ threads has a fixed degree of aggressiveness. ``Conservative'' threads make emphasis on diversification (global search) and are used for increasing the probabilities of finding a feasible solution, even if the parameter space is rugged or weakly structured. ``Aggressive'' threads make emphasis on intensification (local search) and speed up the calculations in smoother areas. Communication, which takes place at fixed time intervals, enables each thread to benefit from the knowledge gathered by the others. Thus this strategy has several degrees of freedom that have to be fixed: the time between communication ($\tau$), the number of threads ($\eta$), and the strategy adopted by each thread. These adjustments should be chosen carefully depending on the particular problem we want to solve. Some guidelines for doing this can be found in \cite{Villaverde_et_al:12} and in the supplementary material accompanying this paper.

	\section*{Implementation}
	MEIGO runs on Windows, Mac, and Linux, and provides implementations in both Matlab and R. So far, MEIGO includes: (i) eSS (Enhanced Scatter Search, \cite{Egea_Marti_Banga_10}), for solving cNLP and MINLP problems, and (ii) VNS (Variable Neighbourhood Search), following the implementation described in \cite{Hansen_Mladenovic_Moreno-Perez:10}, to solve IP problems (see Figure 1). Cooperative parallel versions (CeSS, CVNS), which can run on multicore PCs or clusters, are also included. Cooperation enhances the efficiency of the methods, not only from the computational time point of view, but also because the threads running in parallel are completely independent they can be customized to cover a wide range of search options, from aggressive to robust ones. In a sense, the cooperation as it has been designed, acts as a combination of different metaheuristics since each of the threads may present a different search profile.
	Four different kernel functions per method are included depending on the programming language chosen and the parallelization capabilities. Parallel computation in Matlab is carried out making use of the jpar tool \cite{Karbowski_et_al:08}. Parallel computation in R can be performed using the package snowfall \cite{Knaus:10}. 
	
\begin{figure}[!ht]
    \centering
        \includegraphics[width=0.75\textwidth]{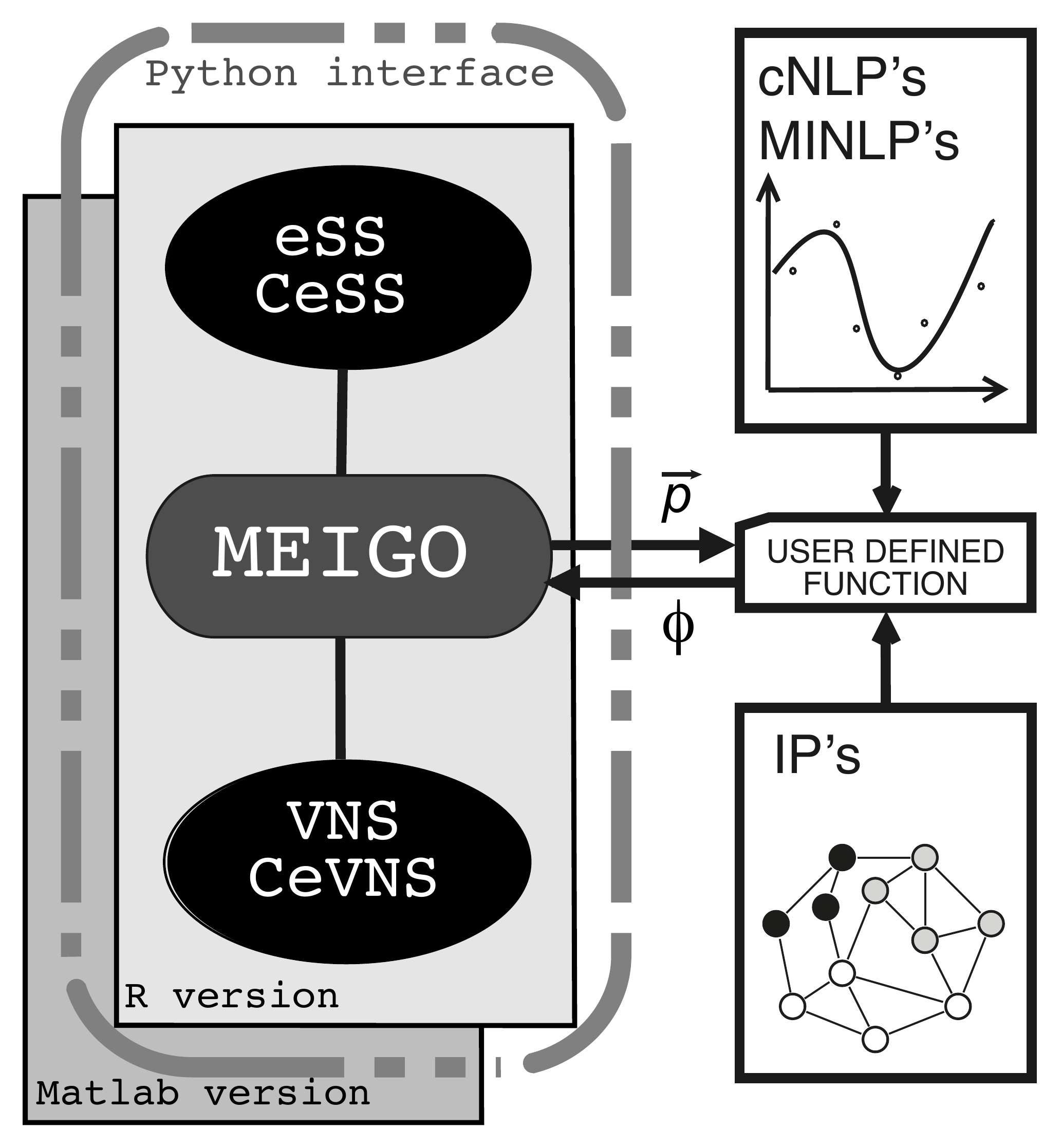}
    \caption{Figure depicting the global structure of MEIGO}
\end{figure}

The methods implemented in MEIGO consider the objective functions to be optimized as black-boxes, with no requirements with respect to their structure. The user must provide a function that can be externally called for evaluation, accepting as input the variables to be estimated, and providing as output the objective function value, $\phi$, associated to the input parameters and the values of the constraints for constrained problems. For either method (eSS and VNS), the user must define a set of compulsory fields (e.g., the name of the objective function, the bounds in the parameters, the maximum number of function evaluations) and a set of options to explore others than the default ones. After each optimization, all the necessary results are stored in data files for further analysis with the tools provided by the host platforms.

Importantly, MEIGO is an open optimization platform in which other optimization methods can be implemented regardless their nature (e.g., exact, heuristics, probabilistic, single-trajectory, population-based, etc.)

\section*{Illustrative examples}
To illustrate the capabilities of the methods presented here, a set of optimization problems, including cases from systems biology and bioinformatics, have been solved and they are presented as case studies. The examples include (i) a set of state of the art benchmark cases for global optimization (from the Competition on Large Scale Global Optimization, 2012 IEEE World Congress on Computational Intelligence), (ii) a metabolic engineering problem based on a constraint-based model of \textit{E. coli}, and (iii) training of logic models of signaling networks to phospho-proteomic data \cite{Saez-Rodriguez_et_al:09}. The corresponding code for these examples is included in the distribution of the MEIGO software.

\subsection*{Large-Scale Continuous Global Optimization Benchmark}
These are benchmark functions used in the Special Session on Evolutionary Computation for Large Scale Global Optimization, which was part of the 2012
IEEE World Congress on Computational Intelligence (CEC@WCCI-2012). These objective functions can be regarded as state-of-the-art benchmark 
functions to test numerical methods for large-scale (continuous) optimization. Information about the functions as well as computer codes can be downloaded 
from \url{http://staff.ustc.edu.cn/~ketang/cec2012/lsgo_competition.htm}.
Some of these functions were previously solved by \cite{Villaverde_et_al:12} using \textit{CeSS}, a cooperative version of the Enhanced Scatter Search 
metaheuristic implemented in Matlab and available within \textit{MEIGO}. Large-scale calibration of systems biology models were also presented and solved in that paper.
Here we present the solution of 3 of these functions (i.e., \textbf{f10}, \textbf{f17} and \textbf{f20}) using the R version of \textit{CeSS} used by \textit
{MEIGO}. The convergence curves for the solution of these benchmark functions in R are coherent with those presented by \cite{Villaverde_et_al:12} solved with Matlab, and the results are competitive with the reference results for these functions presented in \url{http://staff.ustc.edu.cn/~ketang/cec2012/lsgo_competition.htm}. The convergence curves corresponding to these results are presented in Figures 2-4.

\begin{figure}[!ht]
    \centering
        \includegraphics[width=0.75\textwidth]{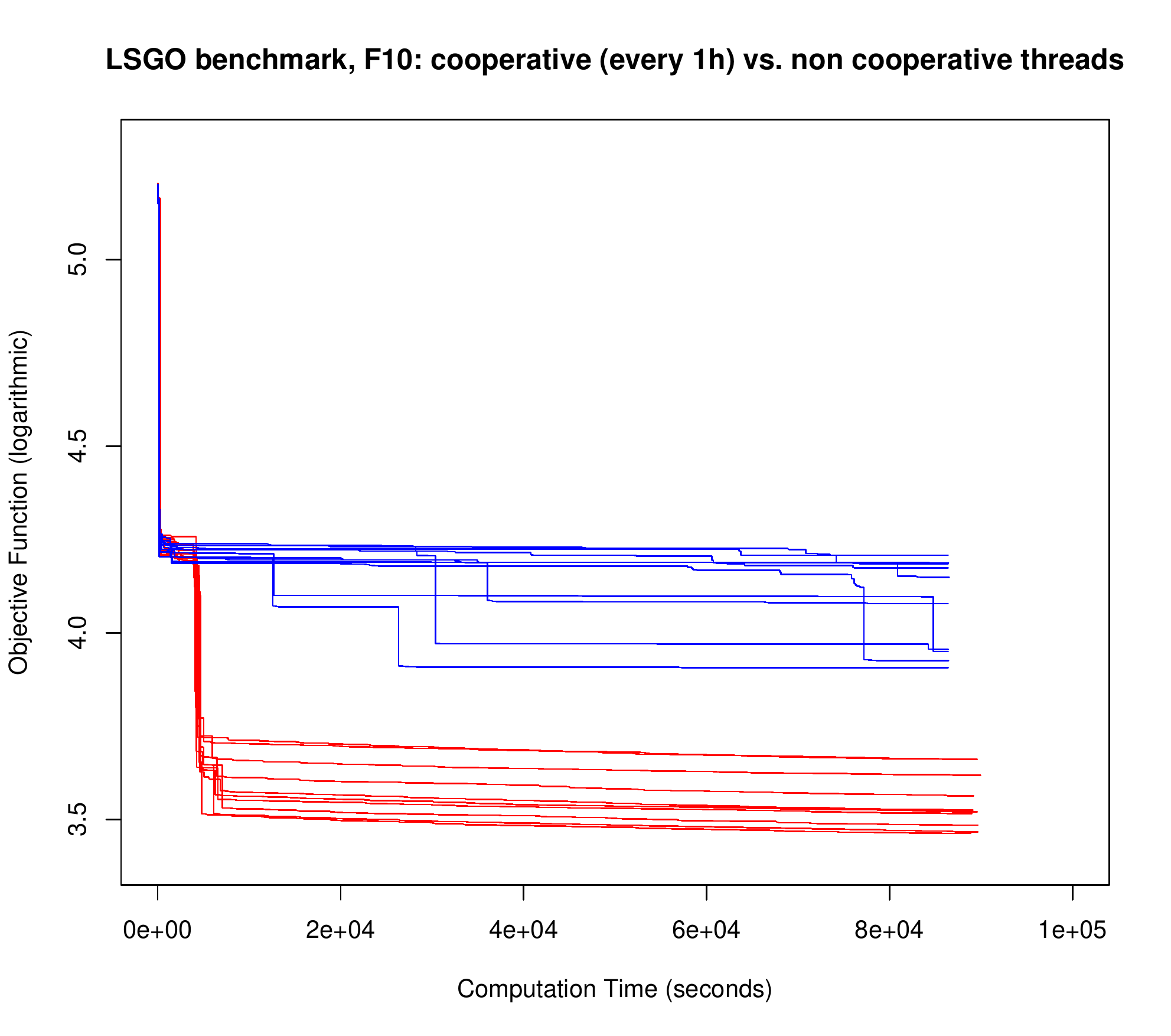}
    \caption{Convergence curves for f10 function}
\end{figure}

\begin{figure}[!ht]
    \centering
        \includegraphics[width=0.75\textwidth]{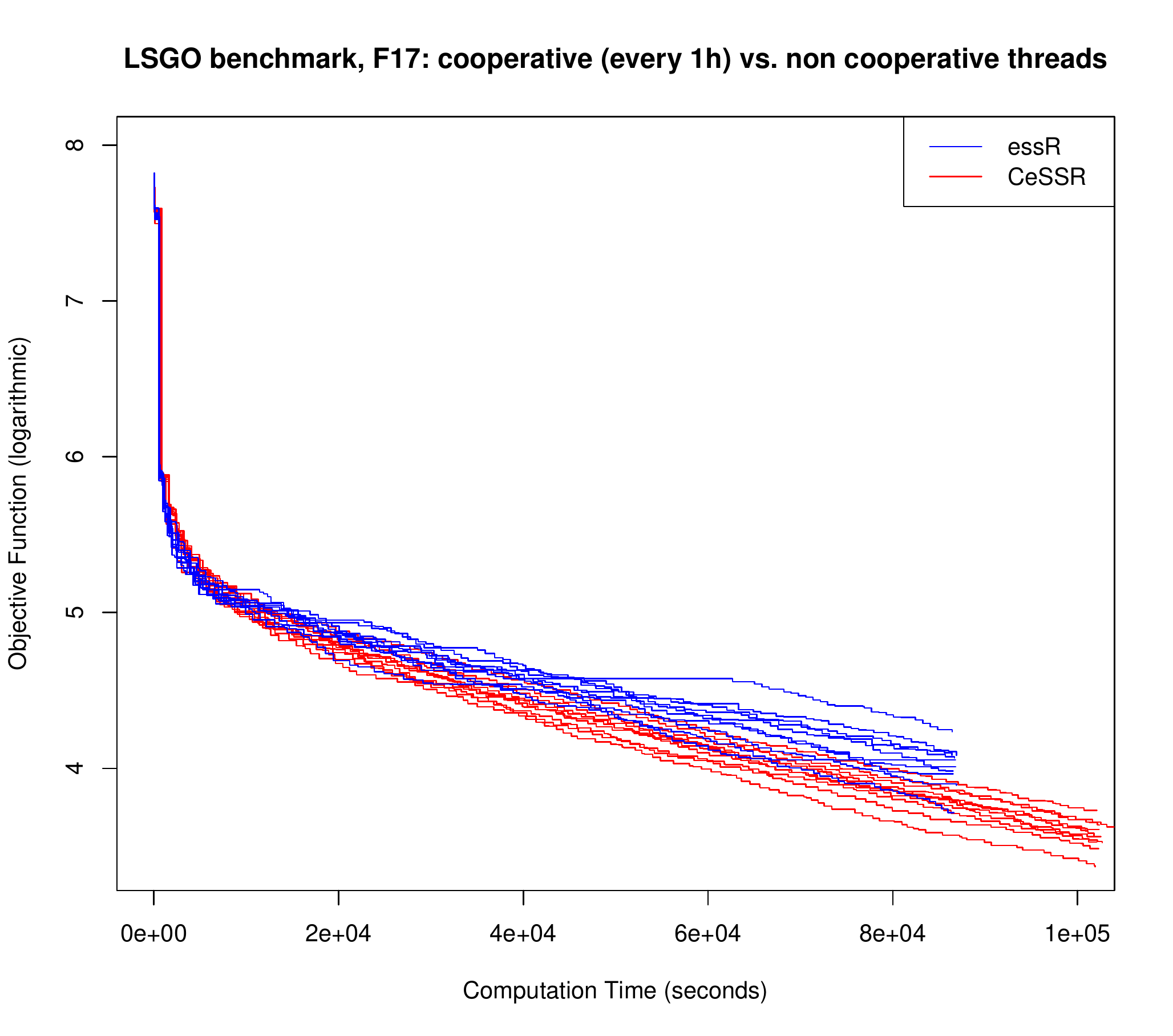}
    \caption{Convergence curves for f17 function}
\end{figure}

\begin{figure}[!ht]
    \centering
        \includegraphics[width=0.75\textwidth]{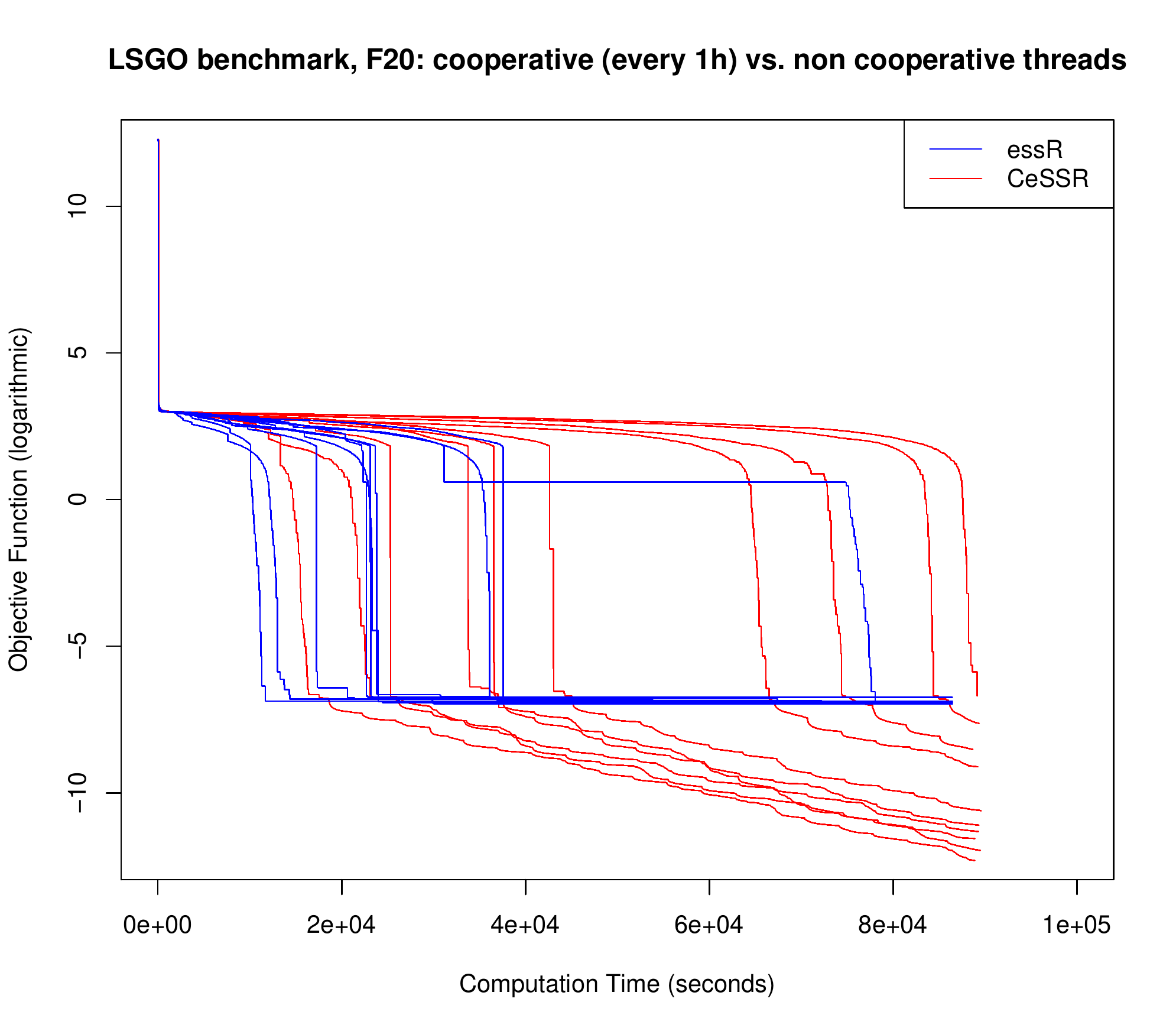}
    \caption{Convergence curves for f20 function}
\end{figure}

\subsection*{Integer optimization benchmark problems}
A set of integer optimization problems arising in process engineering and coded in AMPL (\textit{A Modeling Language
for Mathematical Programming}) were solved using the Matlab version of VNS and making use of the AMPL-Matlab interface files provided by Dr. Sven Leyffer in his web page \url{http://www.mcs.anl.gov/~leyffer/macminlp/}. VNS solved all the problems and, in some cases, achieved a better solution than the best reported one. A summary of the tested problems is presented in Table 1. These benchmarks have been solved using the Matlab version of MEIGO under Windows only since the dynamic library to access AMPL files runs on Windows.

\begin{table}[!ht]
\centering
\begin{tabular}{|c|c|c|c|c|}
\hline
\hline
\multirow{2}{*}{Name} & \multirow{2}{*}{nvar} & \multirow{2}{*}{Ref.} & Best reported  & Best VNS\\
& & & solution & solution\\
\hline
\hline
geartrain & 4 & \cite{Sandgren:90} &7.78e-7 & \textbf{2.70e-12}\\
\hline
mittelman &16 & -  & 13.0 & 13.0\\
\hline
trimlon2 &8 &\multirow{3}{*}{\cite{Harjunkoski_et_al:98}} &5.3 & 5.3\\
\cline{1-2}
\cline{4-5}
trimlon4 &24 & &11.3 & \textbf{8.3}\\
\cline{1-2}
\cline{4-5}
trimlon5 &35 & &12.1 & \textbf{10.6}\\
\hline
\end{tabular}
\end{table}

\subsection*{Metabolic engineering example}
In this section we illustrate the application of the VNS (variable neighbourhood search) algorithm with a metabolic 
engineering problem. Here VNS was used to find a set of potential gene knock-outs that will maximize the production of a given metabolite of interest. The objective function is given by flux-balance analysis (FBA) where a steady-state model is simulated by means of linear programming (LP). FBA assumes that cells have a biological objective often considered as growth rate maximization, minimization of ATP consumption or both.

In this example we considered a small steady-state model from \textit{E. coli} central carbon metabolism available 
at \url{http://gcrg.ucsd.edu/Downloads/EcoliCore}.
Here the metabolite of interest is succinate and we considered the biological objective as biomass maximization.
To solve the inner FBA problem we used openCOBRA (\url{http://opencobra.sourceforge.net/}) with Gurobi as LP solver (
\url{http://www.gurobi.com/}).
As for the problem encoding, 5 integer variables were chosen as decision variables, one for each possible gene knock-out. Each of these variables was allowed to vary from 0 (no knock-out) to 52, the total number of possible genes to be knocked-out. Repeated solutions were filtered by the objective function. 

Additionally we also implemented and solved the problem with a genetic algorithm from the Matlab Global Optimization Toolbox. The point here was to cross-check the VNS results, not to perform an extensive comparison between the performance of GA and VNS. However we found out that for our particular problem and encoding, VNS achieved the optimal solution more often (see Figures 5 and 6). The Wilcoxon rank sum test with continuity correction for comparing means provides a $p-value=0.06753$ (or $0.02104$ if we remove the outlier VNS solution) showing that the solutions provided by VNS are significantly better. Note that the GA was used out of the box (with default settings). Results can vary using other encodings and further tuning of the search parameters. In any case, the purpose was to illustrate how this class of problems can be easily solved using VNS.

\begin{figure}[!ht]
    \centering
        \includegraphics[width=0.75\textwidth]{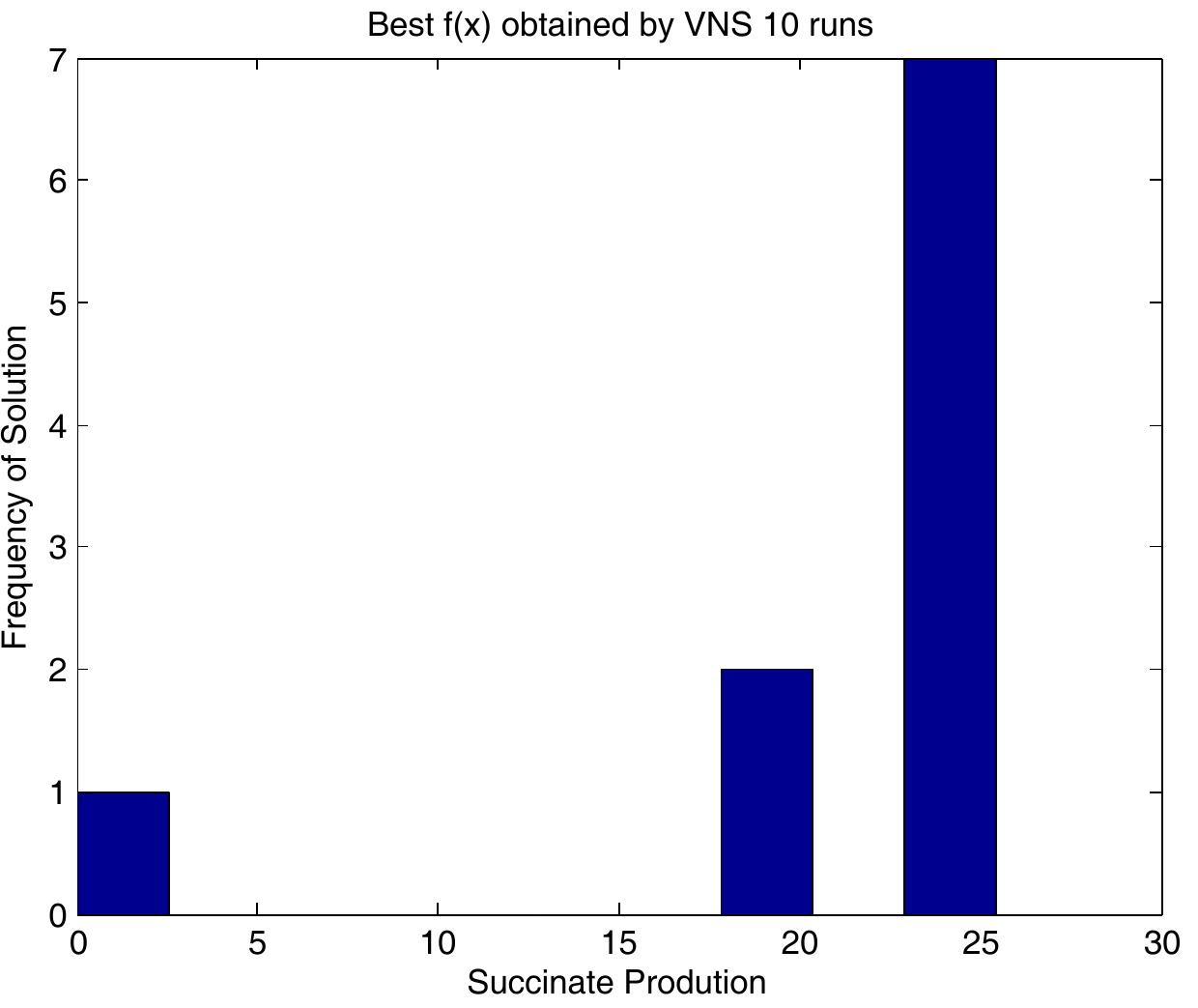}
    \caption{Histogram of the solutions obtained by VNS over 10 runs for the metabolic engineering
    example.}
\end{figure}

\begin{figure}[!ht]
    \centering
        \includegraphics[width=0.75\textwidth]{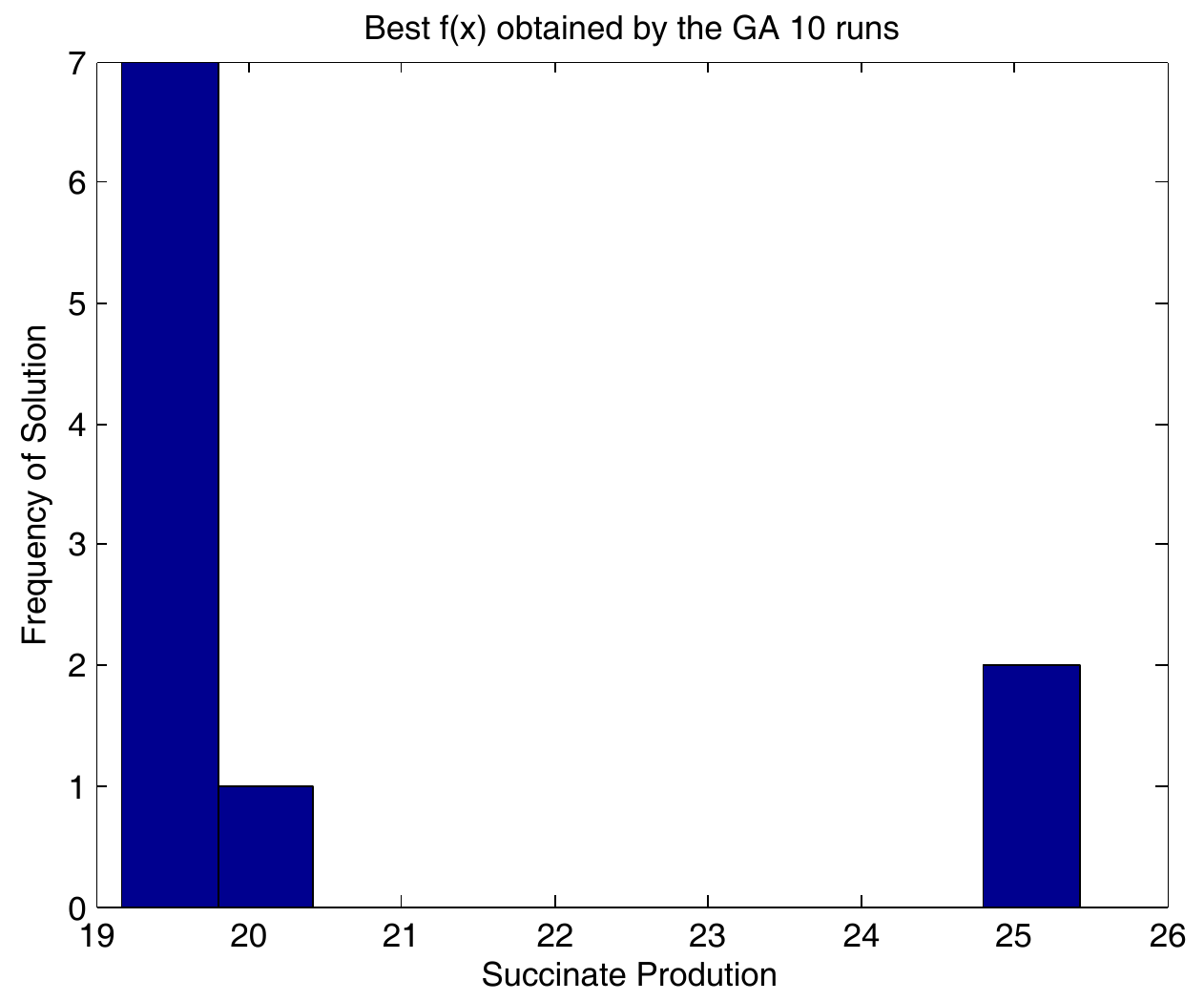}
    \caption{Histogram of the solutions obtained by the genetic algorithm over 10 runs for the metabolic
    engineering example.}
\end{figure}

\subsection*{Training of logic models of signalling networks to phospho-proteomic data}
In this section we compare the performance of variable neighborhood search (VNS) and a discrete genetic algorithm (GA) implementation for the training a logic model of a signalling network to phospho-proteomic data \cite{Saez-Rodriguez_et_al:09}. 

The problem is formulated as follows: starting from a directed graph (a so-called prior knowledge network) containing the interactions and its sign (activating or inhibitory), we can obtain an expanded hypergraph containing all the possible gates, where edges with two or more inputs (hyperedges) represent a logical disjunction (AND gate) and single edges represent a logical conjunction (OR gate). To calibrate such models, the authors formulated the inference problem as a binary multi-objective problem, where the first objective corresponded to how well the model described the experimental data and the second consisted in a complexity penalty to avoid over-fitting:

\begin{equation}
 {\theta(P)=\theta_f+\alpha \cdotp \theta_s(P)}
\end{equation}

\noindent where \( {\alpha(P)} \) is the mean squared error and \( {\alpha \cdot \theta (P)} \) is the product between a tunable parameter $alpha$ and a function denoting the model complexity (AND gates and OR gates receive twice the penalty of a simple activating or inhibiting edge).  

Noticeably, the binary implementation of this problem contains redundant solutions in the search space. This can be addressed by compressing the search space into a reduced set containing only the smallest non-redundant combinations of hyperedges \cite{Saez-Rodriguez_et_al:09} (equivalent to the Sperner hypergraph). By doing this, the problem is transformed from binary to an integer programming problem, that was solved in \cite{Saez-Rodriguez_et_al:09} using a genetic algorithm.

Here, we implemented this benchmark by using the Matlab version of Cell Net Optimizer (CNO or CellNOpt, available at \url{http://www.cellnopt.org/downloads.html}).
The prior-knowledge network and data-set are also publicly available and thoroughly described at \url{http://www.ebi.ac.uk/~cokelaer/cellnopt/data/ExtLiverPCB.html}. 

In order to assess the performance of both algorithms we solved each problem 100 times using VNS and the GA implementation from CNO. In the allowed time budget, VNS has returned solutions that were on average better than those found by the GA (see Figure 7). The Welch Two Sample t-test for comparing means provides a $p-value=3.45e-14$, which clearly shows that VNS outperforms GA for this problem.

\begin{figure}[!ht]
    \centering
        \includegraphics[width=0.75\textwidth]{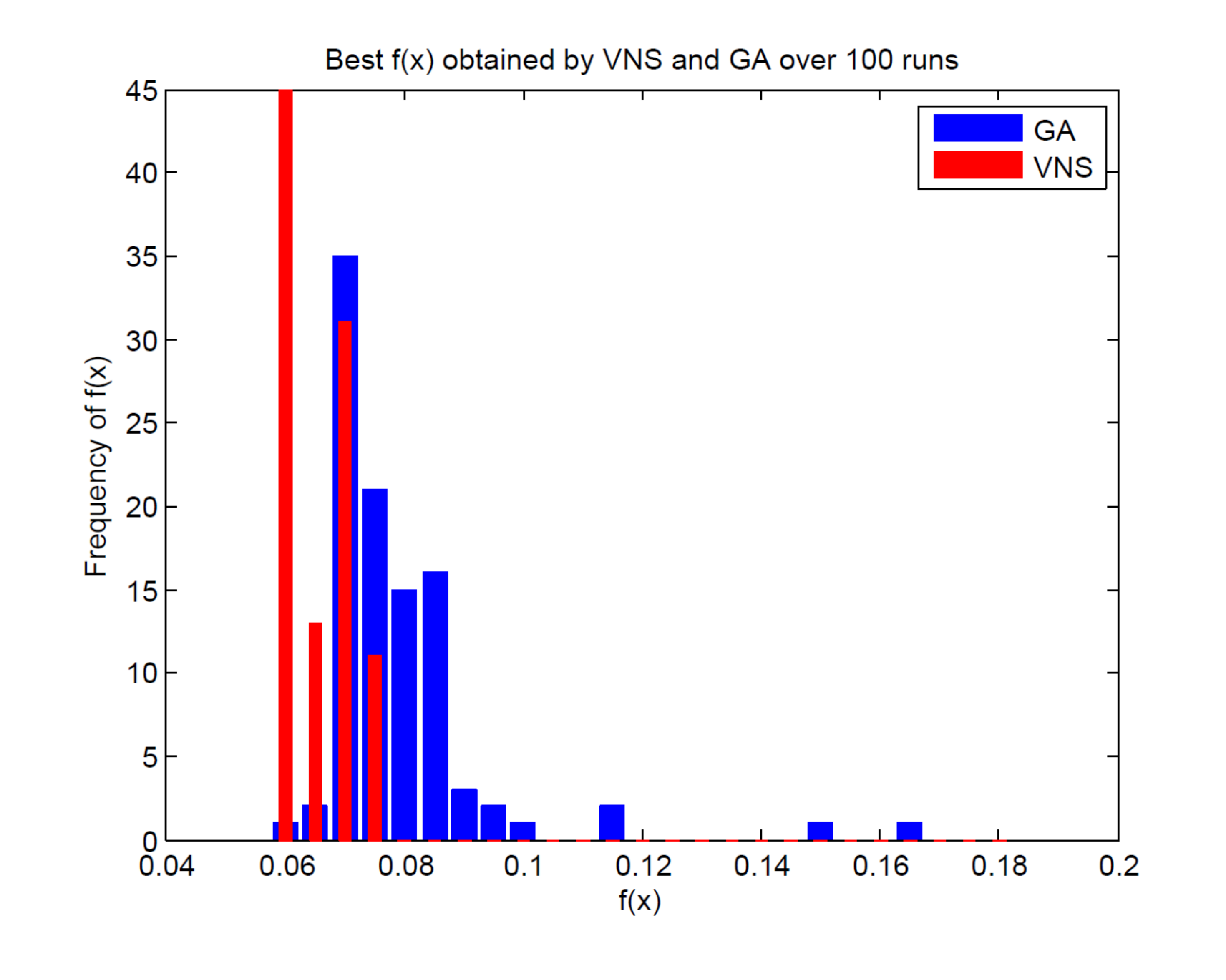}
    \caption{Histogram of solutions obtained by VNS and GA over 100 runs for the logic model example.}
\end{figure}

Since both methods are sensible to the tuning parameters, we tried to tune both algorithms fairly. Also, we note that the solution of this problem in its original, binary
implementation can be solved using deterministic methods based either on Integer Linear Programming \cite{Mitsos_et_al:09,Sharan_Karp:12} or Answer Set Programming \cite{Guziolowski_et_al:13}.

\section*{Conclusions}
Here, we present MEIGO, a free, open-source and flexible package to perform
global optimization in R, Matlab, and Python. It includes advance meta-heuristic
methods, and we are planning to add further methods, such as Bayesian Inference
methods \cite{Wilkinson:07,Eydgahi_et_al:13}. Furthermore, its modular nature (Figure 1), allows the connection to existing optimization methods.

\section*{Availability and requirements}
Project home page: http://www.iim.csic.es/~gingproc/meigo.html\\
Operating system(s): Windows, Linux, Mac OS X\\
Programming language: Matlab 7.5 or higher and R 2.15 or higher\\
Licence: GPLv3

	\bibliography{meigo_arxiv}
	\bibliographystyle{unsrt}

\end{document}